\documentclass[12pt]{amsproc}
\usepackage{amsmath,amsthm}
\usepackage{amsfonts,amssymb}
\newtheorem{thm}{Theorem}[section]

\theoremstyle{remark}

\begin{document}

\title[Characterisation of Fourier transform]
{A Characterisation of the Fourier transform on the
Schwartz-Bruhat space of locally compact Abelian groups}
\author{R. Lakshmi Lavanya}

\address{Department of Mathematics\\ Indian Institute
of Science Education and Research\\Tirupati-517 507}
\email{rlakshmilavanya@iisertirupati.ac.in}

\date{\today}


\begin{abstract}
We obtain a characterisation of the Fourier transform on the space
of Schwartz-Bruhat functions on locally compact Abelian groups.
The result states that any appropriately additive bijection of the
Schwartz space onto itself, which interchanges convolution and
pointwise products is essentially the Fourier transform. The proof
of this result is very similar to that obtained by the author
recently for the Euclidean Fourier transform.
\end{abstract}

\maketitle

\section{Introduction}
\setcounter{equation}{0}

The properties of the Fourier transform on various locally compact
groups with respect to different operations on function spaces on
these groups have been well understood. The interaction between
the Fourier transform and the translations on the groups, and with
certain products on the functions defined on these groups have
been used to obtain characterisations of the
Fourier transform. For more details, refer to \cite{AAM1},\cite{AAM2},\cite{E},\cite{F},\cite{J}-\cite{LT} and the references therein.\\

In 2008, S. Alesker, S. Artstein-Avidan and V. Milman\cite{AAM2}
characterised the Euclidean Fourier transform with the hypothesis
involving only algebraic properties of the map on the space of
tempered distributions on $\mathbb{R}^n.$

\begin{thm}
Assume that $T:\mathcal{S}(\mathbb{R}^n) \rightarrow
\mathcal{S}(\mathbb{R}^n)$ is a bijection which admits a bijective
extension $T': \mathcal{S}'(\mathbb{R}^n) \rightarrow
\mathcal{S}'(\mathbb{R}^n)$ such that for all $f\in
\mathcal{S}(\mathbb{R}^n)$ and $\varphi \in
\mathcal{S}'(\mathbb{R}^n),$ we have
$$T(f\ast \varphi) = T(f) T(\varphi) \textrm{ \ and \ } T(f\cdot \varphi) = T(f) \ast T(\varphi).$$
Then $T$ is essentially the Fourier transform: that is, for some
matrix $B\in GL(n,\mathbb{R})$ with $|det \ B|=1,$ we have either
$Tf= \mathcal{F}(f\circ B)$ or $Tf=\mathcal{F}(\overline{f\circ
B})$ for all functions $f\in \mathcal{S}(\mathbb{R}^n).$
\end{thm}
As the authors of the above result had remarked, the hypotheses of
this result involves only \textit{algebraic} properties of the map
on the class of tempered distributions, while the conclusion
states
that the map is essentially the Fourier transform.\\

The above result led to a characterisation of the Fourier
transform on the Schwartz space of the Heisenberg group\cite{LT}.
This result did not involve the tempered distributions in its
hypothesis. The anonymous referee of \cite{LT} asked an
interesting question, namely, if for the Euclidean Fourier
transform also, one could obtain a characterisation without any
assumptions on the tempered distributions. This question was
answered in affirmation in \cite{La} via the following result:\\

\begin{thm}\label{LaT} Let $T:\mathcal{S}(\mathbb{R}^n) \rightarrow \mathcal{S}(\mathbb{R}^n)$ be a bijection satisfying the following conditions for all functions
$f,g \in \mathcal{S}(\mathbb{R}^n):$
\begin{description}
\item[(a)]  $T(f+\overline{g}) = T(f)+[T(g)]^*,$ where $[Tg]^*(x)
= \overline{Tg(-x)}, \ x\in ~\mathbb{R}^n.$
 \item[(b)] $T(f\cdot g)=
T(f) \ast T(g )$, \item[(c)] $T(f\ast g)= T(g)\cdot T(g).$
\end{description} Then there exists a matrix $B\in GL(n,\mathbb{R}),$ with $|det \ B|=1$ such that
either $Tf = \mathcal{F}(f\circ B)$ for all $f\in
\mathcal{S}(\mathbb{R}^n),$ or $Tf = \mathcal{F}(\overline{f\circ
B})$ for all $f\in \mathcal{S}(\mathbb{R}^n).$
\end{thm}

%

\noindent Hereon, we write LCA groups for locally compact Abelian
groups.\\

 It is interesting to know that the above result holds
good for a general LCA group, with minor changes to its proof. In
this paper, we state and prove this characterisation of the
Fourier transform on LCA
groups. We will discuss analogous characterisations on non-Abelian groups in a future work.\\

Before stating our results, we recall the Schwartz-Bruhat space of
functions on a locally compact Abelian group.\\

F. Bruhat\cite{Br} extended the notion of a smooth function to a
large class of groups, which encompasses the LCA groups. 1n 1975,
M.S. Osborne characterised the Schwartz-Bruhat space of functions
on a LCA group in terms of the asymptotic behavior of the function
and its Fourier transform.\\


 Let $G$ be a locally compact Abelian group, and $\widehat{G},$ the unitary dual of
 $G.$ We use lowercase alphabets $x,y$ for elements of $G,$ the letters $f,g,h$ for functions on $G,$ and the lowercase Greek letters $\alpha, \beta$ for complex
 numbers. We denote the product of two elements $x,y \in G$ as $xy.$ We use the notation $dx$ and $d\xi$ for the Haar measures on
 $G $ and $\widehat{G}$ respectively.
 For $x\in G,$ and $\xi \in \widehat{G},$ we denote this duality by $\langle x,\xi \rangle.$ For an integrable function $f$ on $G,$ its Fourier transform is
 defined as
 $$\widehat{f}(\xi) = \int \limits_G f(x) \ \overline{\langle x,\xi \rangle} \ dx,$$
 where $dx$ denotes the Haar measure on $G.$\\

 Let $f^*(x) =\overline{f(x^{-1})}.$ Then the Fourier transform
 satisfies $$\widehat{f^*}(\xi) = [\widehat{f}]^*(\xi), \textrm{ \ for \ } \xi \in
 \widehat{G}.$$

 We now recall the definition of $\mathcal{S} = \mathcal{S}(G),$ the
 Schwartz-Bruhat space of functions on $G.$ Suppose $G$ is an
 elementary group, i.e., $G$ is of the form $G= \mathbb{R}^n \times \mathbb{T}^m \times
 \mathbb{Z}^k \times F,$ where $n,m,k \in \mathbb{N}_0=\{0,1,2,\cdots\},$ and $F$ is a
 finite Abelian group. A function $f: G\rightarrow \mathbb{C}$ is said to
 belong the Schwartz-Bruhat space if $f$ satisfies the following
 conditions:
 \begin{description}
 \item[(a)] $f\in \mathcal{C}^\infty(G).$
 \item[(b)] $P(\partial)f\in L^\infty(G)$ for all polynomial differential
 operators $P(\partial),$ where the polynomial is in $\mathbb{R}^n \times
 \mathbb{Z}^k$
 variables.
 \end{description}
 For any locally compact Abelian group $G,$ we have
 $$\mathcal{S} (G) = \lim \limits_{\rightarrow} \ \mathcal{S} (H/K),$$
 where the direct limit is taken over all pairs $(H,K)$ of subgroups
 of $G$ such that
 \begin{description}
    \item[(i)] The subgroup $H$ is open and compactly generated,
    \item[(ii)] The subgroup $K$ is compact
    \item[(iii)] The quotient $H/K$ is a Lie group.
 \end{description}

 For a more detailed definition of the space $\mathcal{S} (G),$ we refer the
 reader to \cite{Wa}, where the following
 result was proved.
 \begin{thm}\cite{Wa}
 The Fourier transform maps $\mathcal{S} (G)$ isomorphically onto
 $\mathcal{S} (\widehat{G})$ and $\mathcal{S} (G)$ is dense in $L^1(G).$
 \end{thm}

 For a function $f:G \rightarrow \mathbb{C},$ the support of $f,$
denoted $Supp \ f,$ is defined as
$$Supp \ f := Closure({\{x\in G : f(x) \neq 0\}}).$$

 Let $\mathcal{C}_c^\infty  = \mathcal{C}_c^\infty (G): = \{f\in \mathcal{S} (G) : Supp \ f \textrm{ \ is  \ compact}\}.$ Then the
 space $\mathcal{C}_c^\infty $ is dense in $\mathcal{S} (G).$\\

 For functions $f,g:G \rightarrow C,$ let
 \begin{eqnarray*}
   (fg)(x)& =& f\cdot g(x) = f(x) g(x) \\
    (f \ast g)(x) &=& \int\limits_G f(xy^{-1}) \ g(y) \ dy
 \end{eqnarray*}
Then $\mathcal{S} (G)$ forms an algebra with the above products of functions.\\

Also, the Fourier transform satisfies
$$\widehat{(f \ast g)} = \widehat{f} \ \widehat{g}, \textrm{ \ \ and  \ \ } \widehat{(f\cdot g)} = \widehat{f} \ast \widehat{g} , \textrm{ \ for  \ all \ } f,g \in \mathcal{S} (G).$$

We now proceed to our main result.\\
\section{A Characterisation of Fourier transform on $\mathcal{S}(G)$}

Our results are influenced by the those
of Alesker et al.\cite{AAM2} and their interesting proofs.\\

The proof of the following result is very much similar to that of
Theorem \ref{LaT} given in \cite{La}, except for some minor
variations. For the sake of completeness, we present the proof
here.

\noindent Our main result is the following:

\begin{thm} Let $T:\mathcal{S}(G) \rightarrow \mathcal{S}(\widehat{G})$ be a bijection satisfying the following conditions for all functions
$f,g \in \mathcal{S}(G):$
\begin{description}
\item[(a)]  $T(f+g^*) = T(f)+[T(g)]^*,$
 \item[(b)] $T(f\cdot g)=
T(f) \ast T(g )$, \item[(c)] $T(f\ast g)= T(g)\cdot T(g).$
\end{description} Then there exists a measure-preserving homeomorphism $\psi$ on $G$ such that
either $Tf = \widehat{(f\circ \psi)}$ for all $f\in
\mathcal{S}(G),$ or $Tf = \widehat{(\overline{f\circ \psi)}}$ for
all $f\in \mathcal{S}(G).$
\end{thm}

\begin{proof}
For $f\in \mathcal{S}(G),$ we have $Tf \in
\mathcal{S}(\widehat{G}).$ Since the Fourier transform is a
bijection of $\mathcal{S}(G)$ onto $\mathcal{S} (\widehat{G}),$
there exists unique $g\in \mathcal{S}(G)$ with $Tf=\widehat{g}.$
Define a map $U: \mathcal{S}(G)\rightarrow \mathcal{S}(G)$ as
$Uf:=g$ if $Tf=\widehat{g}.$ Then $Tf = \widehat{(Uf)}$ for all
$f\in \mathcal{S}.$ The map $U$ is a bijection of $\mathcal{S}$
onto itself and satisfies the following conditions for all
functions $f,g\in\mathcal{S}(G):$
\begin{enumerate}
\item  $U(f+g^*) = U(f)+[U(g)]^*,$ \item $U(f\cdot g)= U(f) \cdot
U(g )$, \item  $U(f\ast g)= U(f)\ast U(g).$
\end{enumerate}
The theorem is a then a consequence of the following result, which
gives a precise description of the map $U.$
\end{proof}
\begin{thm}\label{LaU} Let $U:\mathcal{S}(G) \rightarrow \mathcal{S}(G)$ be a bijection satisfying the following conditions for all functions
$f,g \in \mathcal{S}(G):$
\begin{enumerate} \item $U(f+g^*) = U(f)+[U(g)]^*,$ \item $U(f\cdot g)=
U(f) \cdot U(g)$, \item $U(f\ast g)= U(f)\ast U(g).$
\end{enumerate} Then there exists a measure-preserving homeomorphism $\psi: G \rightarrow G$ such that
either $Uf = f\circ \psi$ for all $f\in \mathcal{S}(G),$ or $Uf =
\overline{f\circ \psi}$ for all $f\in \mathcal{S}(G).$
\end{thm}

\begin{proof} We prove the result in 12 steps.

\noindent For $x_0\in G,$ define $$C(x_0): = \{f \in
\mathcal{S}(G): x_0\in Supp\ f\}.$$\\
\vspace{-0.8cm}

\noindent \textbf{Step 1.} Let $f,g \in \mathcal{S}.$ If $g=1$ on
$Supp \ f,$ then $Ug=1$ on
$Supp \ Uf.$\\

\noindent \textit{Proof of Step 1.} Since $g=1$ on $Supp \ f,$ we
have $f\cdot g =f.$ This gives $Uf= U(f\cdot g) = Uf \cdot Ug,$
and so $Ug=1$ on the set $\{x : Uf(x) \neq 0\}.$\\

\noindent Let $x\in Supp \ Uf$ with $Uf(x)=0.$ Then there is a
sequence $\{x_k\}_{k\in \mathbb{N}} \subseteq ~G$ with $Uf(x_k)
\neq 0$ for all $k$ and $x_k \rightarrow x$ as $k \rightarrow
\infty.$ Since $Uf(x_k) \neq 0,$ we have $Ug(x_k) =1$ for all $k.$
Hence $Ug(x) = \lim \limits_{k \rightarrow \infty} \ Ug(x_k) = 1.$
Thus $Ug=1$ on $Supp \ Uf.$\\

\noindent \textbf{Step 2.}  If $f \in \mathcal{C}_c^\infty,$ then $Uf \in \mathcal{C}_c^\infty.$\\

\noindent \textit{Proof of Step 2.} Choose $f\in
\mathcal{C}_c^\infty$ such that $f(x_0) \neq 0.$ Choose $g\in
\mathcal{S}$ such that $g=1$ on $Supp \ f.$ By Step 1, $Ug=1$ on
$Supp \ Uf.$
Since $Ug\in \mathcal{S},$ we get $Supp \ f$ is compact.\\

%


\noindent \textbf{Step 3.} For any $x_0 \in G,$ there exists $y_0
\in G$ such
that $Uf\in C(y_0)$ whenever $f \in C(x_0).$\\

\noindent \textit{Proof of Step 3. }Let $E:=\{f\in \mathcal{S}:
f(x_0) \neq 0\}.$ Fix a function $g\in \mathcal{C}_c^\infty$ with
$g(x_0) \neq 0.$ By Step 2, we have $K: = Supp \ Ug$ is compact.\\

\noindent  For $f\in E,$ define $K_f:= K \cap Supp \ Uf.$ For
functions $f_0:=g,f_1,\cdots ,f_k\in E,$ we have $\prod_{j=0}^k \
f_j\not \equiv 0,$ and so $\prod_{j=0}^k \ Uf_j= U(\prod_{j=0}^k \
f_j) \not \equiv 0,$ which gives $\cap _{j=0} ^k \ K_{f_j} \neq
\emptyset.$ This means, the collection $\{K_f: f\in E\}$ of closed
subsets of $K$ has finite intersection property. Since $K$ is
compact, this gives $\bigcap \limits_{f\in E} \ K_f \neq
\emptyset.$ Let $y_0 \in \bigcap
\limits_{f\in E} \ K_f.$\\

\noindent \textbf{Claim.} If $f\in C(x_0),$ then $Uf \in C(y_0).$\\

\noindent \textit{Proof of Claim.} We prove the claim in two
separate
cases.\\

\noindent \textbf{Case 1.} $f(x_0) \neq 0.$\\
 Then
$f$ never vanishes on a neighborhood, say, $V$ of $x_0.$ Let $g\in
\mathcal{S}$ be such that $f\cdot g = 1$ on $V.$ Choose $h\in
\mathcal{S}$ such that $h=1$ on a neighborhood $W$ of $x_0,$ and
satisfies $W\subseteq Supp \ h \subseteq V.$ Since $f\cdot g=1$ on
$Supp \ h,$ by Step 1, we get $U(f\cdot g)=Uf\cdot Ug =1$ on
$Supp\
Uh.$ Since $h\in E,$ by definition, $y_0 \in Supp \ Uh .$ This implies $Uf(y_0) \neq 0,$ and hence $Uf\in C(y_0).$\\

\noindent Note that all our arguments till now can be applied to
the map $U^{-1}$ as well, and so we have proved that $f(x_0)\neq
0$ if and
only if $Uf(y_0) \neq 0.$\\

\noindent A function $f:G\rightarrow \mathbb{C}$ is said to
satisfy the condition $(\star)$ if the following holds:
$$(\star) \ \ \ \ \ \ \ \ \ \ f(x_0)=0 \textrm{ \ if \ and\  only \ if \ } Uf(y_0)=0.\hspace{4cm}$$
By the above discussion, we have that all functions in
$\mathcal{S}$ satisfy condition $(\star).$\\

\noindent \textbf{Case 2.} $f(x_0) =0.$\\

\noindent Suppose $Uf \not \in C(y_0).$ Then there is a
neighbourhood, say $W,$ of $y_0$ such that $Uf$ vanishes
identically on $W.$ Let $h\in \mathcal{S}$ with $Supp \ h
\subseteq W,$ and $h(y_0) \neq 0.$ There exists unique function $g
\in \mathcal{S}$ with $Ug=h.$ Then $U(f\cdot g) = Uf\cdot Ug = Uf
\cdot h \equiv 0.$ This gives $f\cdot g
\equiv 0.$ \\

\noindent On the other hand, since $Ug(y_0) = h(y_0) \neq 0,$ by
Condition $(\star),$ we have $g(x_0) \neq 0,$ and so g is never
zero near $x_0.$ Since $x_0\in Supp \ f,$
 this implies $f\cdot g \not \equiv 0,$ a contradiction. Thus $Uf \in
C(y_0).$\\

%
%

\noindent \textbf{Step 4.} Define a map $\varphi : G \rightarrow
G$ as follows: $\varphi (x)=y$ if $Uf\in C(y)$
whenever $f\in C(x).$ Then the map $\varphi $ is well-defined.\\

\noindent \textit{Proof of Step 4.} Suppose for some $x_0 \in G,$
we have $\varphi (x_0) = y_1$ and $\varphi (x_0) = y_2$ with
$y_1\neq y_2.$ Let $V_1$ and $V_2$ be disjoint neighborhoods of
$y_1$ and $y_2,$ respectively. There exists functions $g_1$ and
$g_2$ in $\mathcal{S}$ which are supported
 in $V_1$ and $V_2,$ respectively, such that $g_1(y_1) \neq 0$ and $g_2(y_2) \neq 0.$ Let $f_1,f_2 \in \mathcal{S}$ with $Uf_1=g_1$ and
 $Uf_2=g_2.$ Then $0\equiv g_1\cdot g_2 =Uf_1 \cdot Uf_2 = U(f_1 \cdot f_2)$ and so $f_1\cdot f_2 \equiv
 0.$\\

 On the other hand, as $g(y_j) = Uf_j(y_j)\neq 0$ for $j=1,2,$ we
 have by
Condition $(\star)$ that $f_j(x_0) \neq
 0$
 for $j=1,2,$ which is in contradiction to $f_1\cdot f_2\equiv 0.$\\

\noindent \textbf{Step 5.} The map $\varphi :G \rightarrow G$ is a bijection.\\

\noindent \textit{Proof of Step 5.} The hypotheses of the theorem
hold good for the map $U^{-1}$ as well. Applying the preceding
steps to the map $U^{-1}$ gives rise to a well-defined function,
say,
$\psi :G \rightarrow G.$ Then $\psi  = \varphi ^{-1}$, proving that $\varphi $ is bijective.\\

\noindent Our observations can be summarised as:
$$\varphi (Supp \ f) = Supp \ Uf,\textrm{ \ for \ all \ } f\in \mathcal{S}.$$


%
%


\noindent \textbf{Step 6.} The map $\varphi $ is a homeomorphism of $G$ onto itself.\\

\noindent \textit{Proof of Step 6.} First we prove that $\varphi$
is continuous. Suppose not. Then there exist $x\in G,$ and
sequence $\{x_k\}$ in $G$ with $x_k
\rightarrow x$ as $k \rightarrow \infty,$ but $\varphi( x_k)$ does not converge to $\varphi (x).$\\

\noindent Let $V$ be a neighborhood of $\varphi (x)$ such that
$\varphi( x_k) \not\in V$ for any $k.$ Let $h\in \mathcal{S}$ with
$Supp \ h\subseteq V,$ and $h[\varphi (x)] = 1.$ Let $g\in
\mathcal{S}$ be such that $Ug=h.$ Then $\varphi (x_k )\not \in
Supp \ Ug$ for any $k,$ and so $x_k \not\in Supp \ g$ for any $k.$
This gives
 $g(x_k) = 0$ for all $k,$ implying $g(x)= 0,$ which is not
possible by Condition $(\star),$ since $Ug[\varphi (x)] =1.$\\

We observe that the above argument holds good when the maps $U$
and $\varphi $ are replaced with $U^{-1}$ and $\varphi ^{-1}$
respectively, yielding that
$\varphi : G \rightarrow G$ is a homeomorphism.\\

%
%

\noindent \textbf{Step 7.} The map $\varphi$  satisfies $\varphi (xy) = \varphi (x) \varphi (y)$ for all $x,y \in G.$\\

\noindent \textit{Proof of Step 7.} Suppose $\varphi (xy) \neq \varphi (x) \varphi (y)$ for some $x,y\in G.$\\

\noindent Then there exist disjoint neighborhoods
$V_{xy},V_{x\diamond y}$ with $\varphi (xy) \in ~V_{xy}$ and
$\varphi (x) \varphi (y) \in V_{x\diamond y }.$ By continuity of
the map $\varphi ,$ this gives rise to a neighborhood $W_{xy,1}$
of $xy$ with $\varphi (W_{xy,1}) \subseteq V_{xy}.$ By continuity
of multiplication in $G,$ we get neighborhoods $W_{x,1}, W_{y,1}$
of $x$ and $y,$ respectively, such that $W_{x,1}W_{y,1} \subseteq
W_{xy,1}.$ Thus
\begin{eqnarray} \vspace{-1cm} \varphi (W_{x,1} W_{y,1}) &\subseteq&
\varphi (W_{xy,1}) \subseteq V_{xy}.\label{eq1}
\end{eqnarray}

On the other hand, by continuity of multiplication in $G,$ we have
$\varphi( x) \varphi (y) \in V_{x\diamond y }$ gives neighborhoods
$V_{x,2},V_{y,2}$ such that
\begin{eqnarray} \varphi (x)\in V_{x,2}, \ \varphi (y) \in V_{y,2} \textrm{ \ and \ }
V_{x,2}\ V_{y,2}\subseteq V_{x\diamond y }.\label{eq2}
\end{eqnarray}
 This
implies there exist neighborhoods $W_{x,2}, W_{y,2}$ of $x$ and
$y,$ respectively, with $\varphi (W_{x,2})\subseteq V_{x,2}$ and
$\varphi (W_{y,2}) \subseteq V_{y,2}.$\\

Define $W_x= W_{x,1}\cap W_{x,2}, \ W_y = W_{y,1}\cap W_{y,2}.$
Then
\begin{eqnarray}
  \varphi (W_x)& \subseteq & \varphi (W_{x,2}) \subseteq V_{x,2} =V_x \ (say) \nonumber \\
   \varphi (W_y)& \subseteq & \varphi (W_{y,2}) \subseteq V_{y,2} =V_y \ (say)\nonumber \\
   \varphi (W_x W_y)& \subseteq & \varphi (W_{x,1} W_{y,1})  \subseteq \varphi (W_{xy,1}) \subseteq V_{xy} \label{eq3}
\end{eqnarray}
Choose $f_x,f_y \in \mathcal{S}$ such that $Supp \ f_x \subseteq
W_x, Supp\ f_y \subseteq W_y$ and $f_x \ast ~f_y \not \equiv~0.$
Let $g_x = Uf_x$ and $g_y = Uf_y.$ Then $U(f_x\ast f_y) = g_x \ast
g_y \not\equiv 0.$\\

\noindent We have
\begin{eqnarray}
 \nonumber Supp(g_x \ast g_y) &=&Supp \ U(f_x \ast f_y)  \subseteq  (Supp \ Uf_x) \ (Supp \ Uf_y) \\
\nonumber &= & \varphi (Supp \ f_x) \ \varphi (Supp \ f_y)
\subseteq \varphi( W_x) \ \varphi (W_y)
\\ & \subseteq & V_{x,2}\ V_{y,2} \subseteq V_{x\diamond y } \textrm{   \ \ (by \
\ref{eq2})} \label{eq4}
    \end{eqnarray}
But $Supp \ (f_x  \ast f_y) \subseteq (Supp \ f_x) \ (Supp \ f_y)
\subseteq W_x\ W_y.$ By $(\ref{eq3}),$ this gives
 \begin{eqnarray}
 \nonumber Supp(g_x \ast g_y) &= &Supp (Uf_x \ast Uf_y) = Supp \ U(f_x \ast f_y)\\
  & =& \varphi (Supp (f_x  \ast f_y)) \subseteq  \varphi (W_x \ W_y)\subseteq
 V_{xy}. \label{eq5}
 \end{eqnarray}
From $(\ref{eq4})$ and $(\ref{eq5})$, we get $$Supp(g_x \ast g_y)
 \subseteq V_{x\diamond y } \cap V_{xy} = \emptyset.$$ This gives $g_x\ast
g_y \equiv 0,$ a contradiction.
This proves the multiplicativity of the map $\varphi .$\\

Thus the map $\varphi $ is a homeomorphism of G onto itself, whose
inverse we denote by $\psi.$\\


%
%

\noindent \textbf{Step 9.} 'Extension' of the map $U$ to scalars.\\

\noindent \textit{Illustration of Step 9.} For $f,g \in
\mathcal{S},$ and $\alpha (\neq 0)\in \mathbb{C},$ we have
$$ U(\alpha f)(x)\ Ug(x) = U(\alpha fg)(x) = U(f)(x) \ U(\alpha g)(x) , \ x\in G.$$
Let $h \in \mathcal{S}$ be such that $Uh(x) \neq 0$ for any $x\in
G.$ Then we have
\begin{eqnarray*}
 U(\alpha f)(x)  &=& \frac{U(\alpha h)(x)}{Uh(x)} \ {Uf}(x) \textrm{ \ for \ all \ } f\in \mathcal{S}\\
   &=& m(\alpha ,x) \ Uf(x) \textrm{ \ (say)}.
\end{eqnarray*}
Thus $U(\alpha f)(x) = m(\alpha ,x)\ Uf(x),$ for all $x\in G.$ By
definition, the function $m(\cdot,\cdot)$ is continuous in the
second variable as a function
of $x\in G.$\\

\noindent \textbf{Claim.} The function $m(\cdot,\cdot)$ is
independent of the second
variable.\\
\noindent \textit{Proof of Claim.} For $f,g \in \mathcal{S}, \
\alpha \in \mathbb{C},$ and $x\in G,$ we have
\begin{eqnarray*}
 U(\alpha f\ast g) (x) &=& U(f\ast \alpha g) (x)  \\
 (U(\alpha f)\ast Ug)(x) &=& (Uf \ast U(\alpha g))(x) \\
   \int\limits_{G} m(\alpha ,xy^{-1} )\ Uf(xy^{-1} ) \ Ug(y) \ dy&=&  \int\limits_{G} Uf(xy^{-1} ) \ m(\alpha ,y) \ Ug(y) \ dy \\
\end{eqnarray*}
As the above equation holds good for all functions $Uf,Ug \in
\mathcal{S},$ and $U:\mathcal{S} \rightarrow \mathcal{S}$ is a
bijection, we have for all $f,g \in \mathcal{S},$
   $$\int\limits_{G} [m(\alpha ,xy^{-1} )-m(\alpha ,y)] \ f(xy^{-1} ) \ g(y) \ dy =0, \ x\in G.$$
   Fix $x\in G.$ Let $V$ be a compact symmetric neighborhood of the
   identity $e,$ of $G.$ Choose a function $g\in \mathcal{S}$ such that
   $g=1$ on $V.$ Then for all functions $f\in \mathcal{S}$ with $Supp \ f
   \subset xV,$ we have
\begin{eqnarray*}
  \int\limits_{V} [m(\alpha ,xy^{-1} )-m(\alpha ,y)] \ f(xy^{-1} ) \ dy &=&0. \\
  \textrm{Thus \hspace{2.5cm} } m(\alpha ,xy^{-1} ) -
  m(\alpha ,y) &=& 0 \textrm{\ for all \ } y \in V.
  \end{eqnarray*}
by the continuity of the map $m(\cdot,\cdot)$ in the second
variable. This gives in particular, $m(\alpha ,x) =m(\alpha ,e).$
As $x$ was arbitrary, the above argument gives that the function
$m(\alpha ,x)$ is independent of the second variable $x\in G.$ We
define
$$m(\alpha ): = m(\alpha ,e).$$

\noindent \textbf{Step 10.}  The map $m :\mathbb{C} \rightarrow
\mathbb{C}$ is an additive and multiplicative bijection, which
maps $\mathbb{R}$ onto $\mathbb{R},$ and hence we have either
$m(\alpha) = \alpha$ for all $\alpha\in \mathbb{C},$ or $m(\alpha)
= \overline{\alpha}$ for all
$\alpha\in \mathbb{C}.$\\

\noindent \textit{Proof of Step 10.} Choose $g\in \mathcal{S}, \
\alpha,\beta\in \mathbb{C}$ with $g(x) \neq 0$ for any $x\in G.$
Then by the condition
$(\star),$ we get $Ug[(\varphi(x))]  \neq 0,$ for any $x\in G.$\\

Suppose $m(\alpha) =m(\beta)$ for some $\alpha,\beta\in
\mathbb{C}.$ Then
$$U(\alpha g) (x) = m(\alpha) \ Ug(x)= m(\beta) \ Ug(x) =
U(\beta g)(x) , \ x\in G.$$ Since $U$ is a bijection, this gives $\alpha=\beta.$\\

\noindent By hypothesis(1), we have
\begin{eqnarray*}
 m(\alpha +\overline{\beta})\ Ug(e)  &=& U((\alpha +\overline{\beta})g)(e)
=U(\alpha g+\overline{\beta}g)(e) \\
   &=& U(\alpha g)(e) + \overline{U(\beta\overline{g})}(e) = (m(\alpha ) +m(\overline{\beta} )) \ Ug(e).
\end{eqnarray*}
Since Ug is never zero, we get $m(\alpha +\overline{\beta}) =
m(\alpha )+\overline{m(\beta)}.$ In particular,
$m(\overline{\alpha })
=\overline{m(\alpha )}$ for all $\alpha \in \mathbb{C}.$\\

\noindent Now, hypothesis(2) gives
$$m(\alpha \beta) Ug(e) = U(\alpha \beta g)(e) = m(\alpha )U(\beta g)(e) = m(\alpha ) m(\beta) Ug(e).$$
Again, since $Ug$ is nowhere vanishing, we get $m(\alpha
\beta)=m(\alpha )m(\beta)$
for all $\alpha ,\beta\in \mathbb{C}.$\\

%


%

\noindent \textbf{Step 11.} For $f\in \mathcal{S},$ and $x_0 \in
G,$ we have
$Uf[\varphi (x_0)] = m(f(x_0)).$\\

\noindent \textit{Proof of Step 11.} As before, choose $g \in
\mathcal{S}$ such that $g(x) \neq 0$
for any $x \in G.$ Then $Ug(y) \neq 0$ for any $y \in G.$ \\

\noindent Define $$h(x): = f(x_0) \ g(x) - f(x) \ g(x), \ x\in
G.$$ Then $h\in \mathcal{S}$ and $h(x_0) = 0.$ By Condition
$(\star),$ we have $Uh[\varphi (x_0)] =0.$ This gives
\begin{eqnarray*}
0=Uh[\varphi (x_0)] &=& U[f(x_0) \cdot g - f \cdot g][\varphi (x_0)]\\
 &=& m(f(x_0)) \
Ug[\varphi (x_0)] - Uf[\varphi (x_0)] \ Ug[\varphi (x_0)].
\end{eqnarray*}
 Since $Ug$ is never zero, this gives
$Uf[\varphi (x_0)] = m(f(x_0)).$\\

\noindent Since $\psi =\varphi ^{-1},$ using Step 10, we get that
either
$Uf(x_0) = f[\psi (x_0)]$ or $Uf(x_0) = \overline{f[\psi (x_0)]}.$\\

\noindent Thus we get that the map $U$ is as claimed by our
theorem. It remains to show that $\psi$ is measure-preserving.\\

%

\noindent \textbf{Step 12.} The map $\psi: G\rightarrow G$ preserves measures of subsets of $G.$\\

\noindent \textit{Proof of Step 12.} We have
\begin{eqnarray*}
 (f\ast g) [\psi (x)] = U(f\ast g)(x)  &=& (Uf \ast Ug)(x) \\
   &=& \int\limits_{G} Uf(xy^{-1} ) \ Ug(y) \ dy \\
   &=& \int\limits_{G} f[\psi (xy^{-1} )] \ g[\psi (y)] \ dy\\
      &=& \int\limits_{G} f[\psi (x) \psi(y)^{-1} ] \ g[\psi (y)] \
      dy\\
         &=& \int\limits_{G} f[\psi (x) y^{-1} ] \ g(y) \ d(\varphi(y))
\end{eqnarray*}
This gives $$\int\limits_{G} f[\psi (x) y^{-1} ] \ g(y) \ dy =
\int\limits_{G} f[\psi (x) y^{-1} ] \ g(y) \ d(\varphi(y)),$$ for
all functions $f,g \in \mathcal{S} (G).$ Let $V$ be a compact
symmetric neighborhood of the identity $e$ of $G.$ Choose a
function $g\in \mathcal{S} $ such that $g=1$ on $V.$ Then for all
functions $f\in \mathcal{S} $ which are supported in $xV,$ we get
$$\int\limits_{V} f[\psi (x) y^{-1} )]  \
dy = \int\limits_{V} f[\psi (x) y^{-1} ]  \ d(\varphi(y)).$$ This
implies the map $\varphi$ is preserves the measures of subsets of
$V.$ As $V$ was an arbitrary compact symmetric neighborhood of the
identity, we get the $\varphi,$ and hence $\psi$ is
measure-preserving on subsets of $G.$
\end{proof}

%

\end{document}